\documentclass[twocolumn]{amsart}%
\usepackage{amsmath}
\usepackage{amsfonts}
\usepackage{amssymb}
\usepackage{graphicx}%
\usepackage{wrapfig}
\usepackage[letterpaper, margin=1in]{geometry}
\usepackage{cuted}
\begin{document}

\title{A bug's eye view: The Riemannian exponential map on polyhedral surfaces}
\author{David Glickenstein}
\begin{strip}
  \vspace*{\dimexpr-\baselineskip-\stripsep\relax}
  \centering
  \maketitle
  \vskip\baselineskip
\noindent\makebox[\textwidth]{\rule{1.1\paperwidth}{0.4pt}}
  \vskip\baselineskip
\end{strip}

\begin{quotation}
When a blind beetle crawls over the surface of a curved branch, it doesn't
notice that the track it has covered is indeed curved. I was lucky enough to
notice what the beetle didn't notice.

Albert Einstein - letter to his son Eduard (1922) \cite{einstein}\bigskip
\end{quotation}

\section{Introduction}

What would happen if Einstein's beetle saw in the same way he was able to move? 
If light traveled along a surface to the eye, we would not perceive things in quite the same way.
We will attempt to conceive of what the world would look like
if light traveled in straight lines on the surface of a polyhedron. We will become flatlanders, 
as in E. A. Abbott's famous book \cite{abbott}, except that our Flatland has pockets of curvature 
that come from the failure of vertices of polyhedra to be flat.

Our assumption is that in every direction, we can look out along straight paths (geodesics) and 
see everything that reaches our eye along that path unless it is blocked by something in front of it.
In this way, when we look out, we will see a vast plane of view in every direction and as far as 
the light will carry.

In Riemannian geometry, this view is a representation of the \emph{exponential mapping}. 
The basics of the exponential map can be found in just about any text on Riemannian
geometry, for instance \cite{lee}.

\section{Methods}

In trying to understand how light travels along a surface, it is clear that within a (Euclidean) polygon,
 light travels along straight lines. If a bug is on a polygon, 
 it will see the polygon in the usual (flatland) way. The interesting thing is what happens as the bug looks beyond
 the polygon. Should it look in the direction of an edge, the geometry is such that it will look as though part of the next face (polygon) is folded up into the plane of the first face. The ray from the bug's eye through the two polygons is all visible. The process of folding part of a polygon along an edge into a flat position is called developing the polygon.
 
 But now what happens as the direction of gaze moves and eventually hits a vertex? As the gaze 
 sweeps over a vertex, we can continue to develop polygons into the plane, but any ray that goes through 
 a vertex will cause a cut in the polygon, and stop the development. The rays going through vertices 
 are unpredictable; we will discuss what is really happening later.
 
 We can thus consider the view by starting on the polygon the bug is on, and then recursively developing
  neighboring polygons, cutting along rays that intersect vertices. This produces a ``map view.''
  
  How can we really get a feel for what a bug sees, using our own intuition of what it looks like 
  to see things in a distance? We decided to turn these two-dimensional surfaces (made of polygons
   glued together) into a three-dimensional space that has a very small height. In this way, we can
    consider the bug's eye as slightly above the polygon, and then view what it looks like to be
     the bug. We call this the ``first person view.''

\section{Phenomena}

In this section we show some interesting phenomena that occur on the map and
first person view.

\subsection{Repetition}
\begin{figure}[htbp]
\centering
\begin{tabular}{c}
\includegraphics[width=0.48\textwidth]
{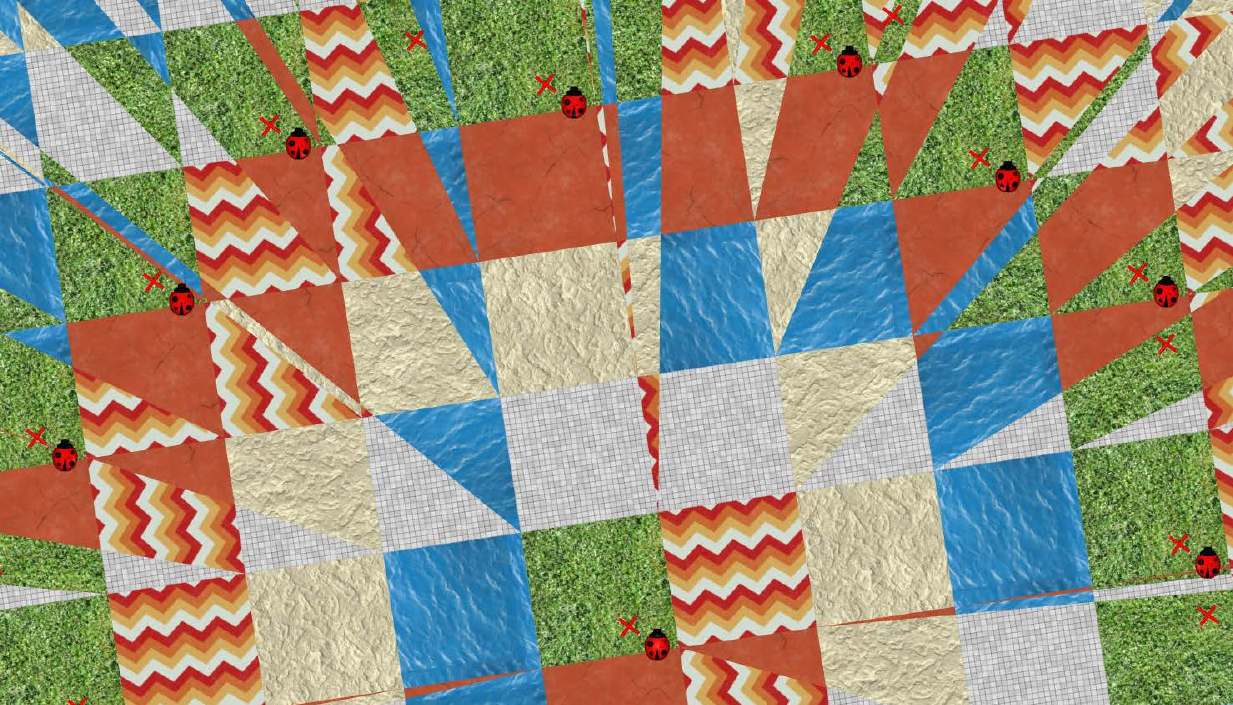}%
\\
{\includegraphics[
width=.48\textwidth
]%
{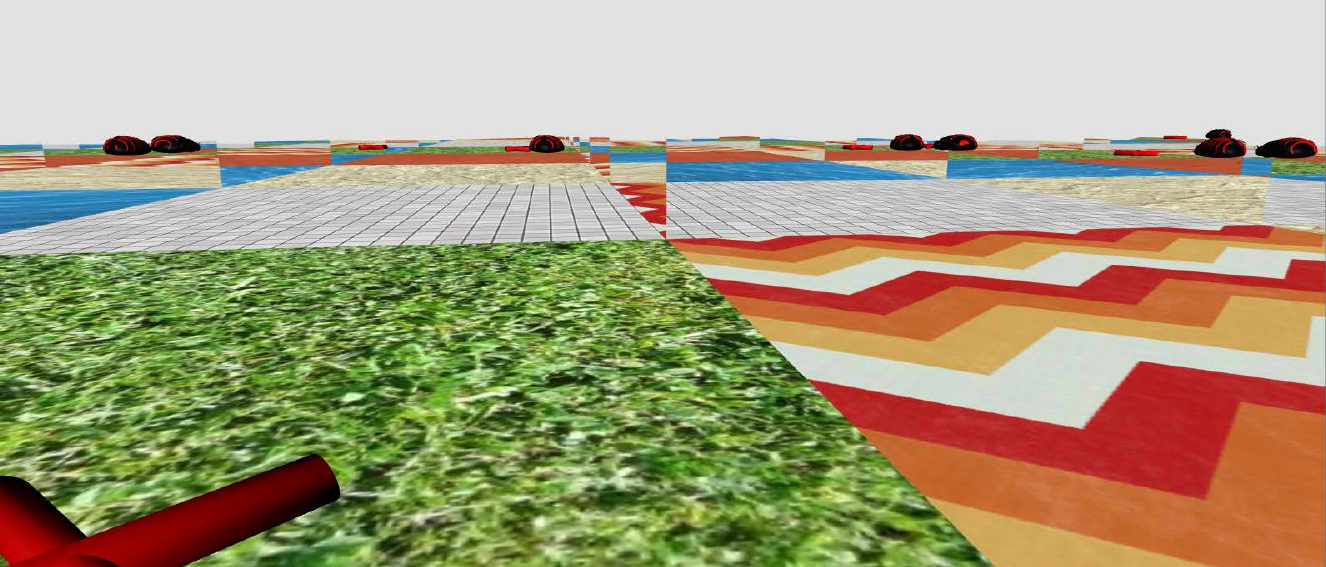}%
}
\end{tabular}
\caption{Repetition on a cube, map and first person views}
\label{fig:mapfirst}
\end{figure}

If one traces a long geodesic, it may curve back around to the same point or
nearby the same point. This phenomenon manifests as seeing multiple images of
the same object along a straight path. The ladybug sees itself multiple times
if it can see far enough relative to the size of the surface. This can occur
on any surface and is the result of short geodesic loops (or actually 1-gons)
on the surface. Figure \ref{fig:mapfirst} shows a picture of repetition. \bigskip%

\subsection{Lensing/positive curvature}
\begin{figure}[htbp]
\centering
\begin{tabular}{c}
\includegraphics[width=0.48\textwidth]
{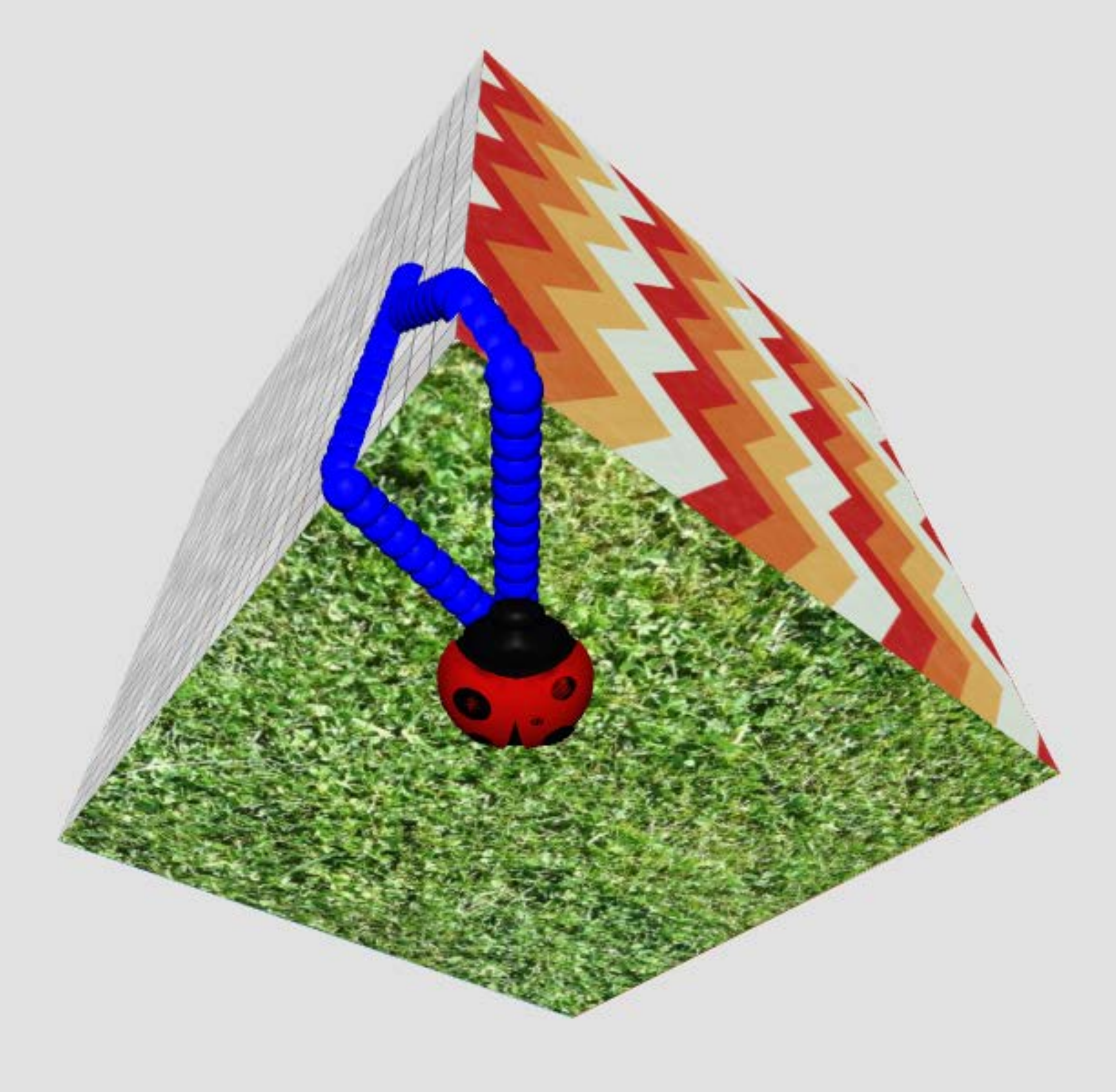}%
\\
\includegraphics[width=0.48\textwidth]
{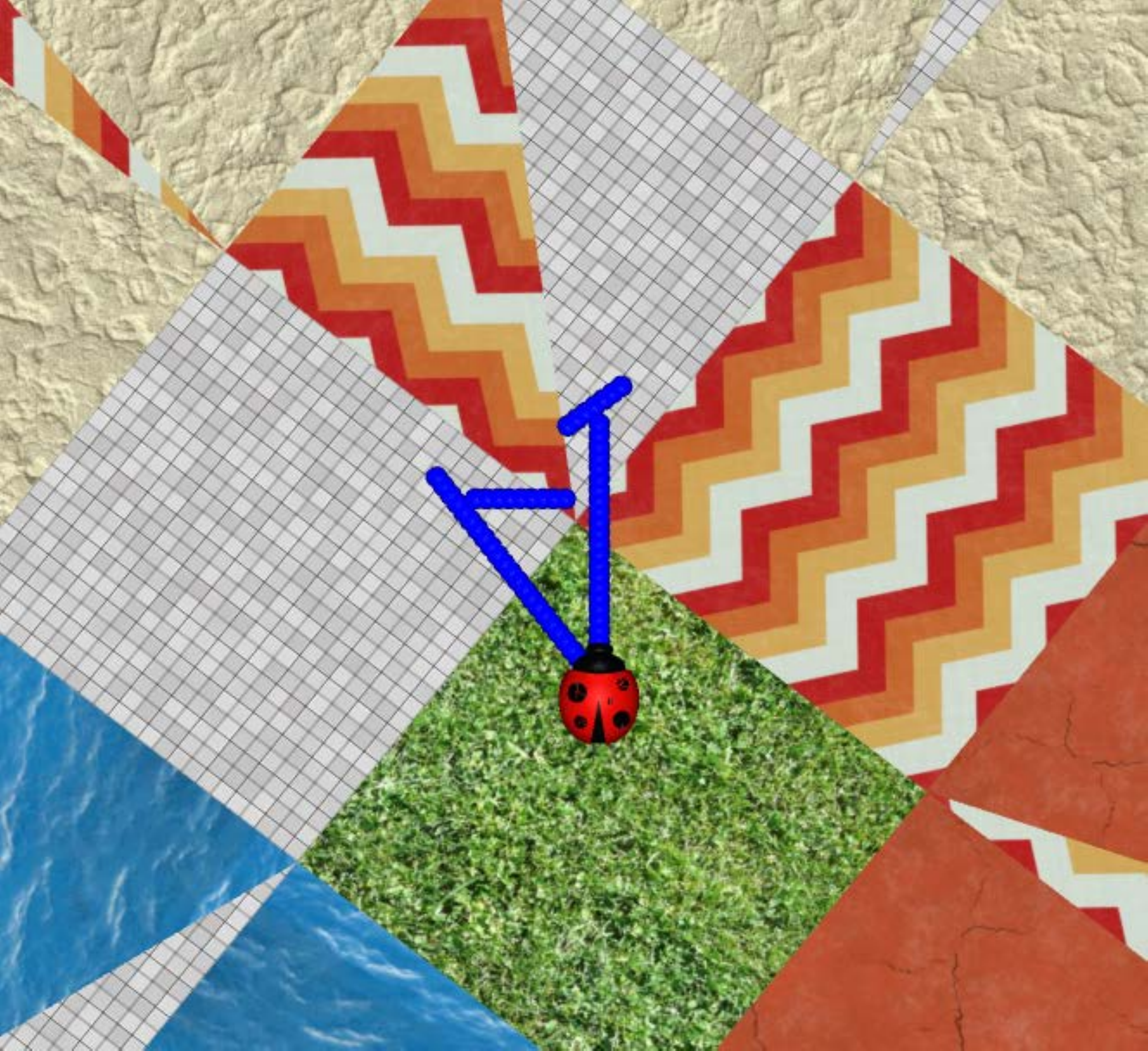}
\end{tabular}
\caption{Geodesics on the cube that curve around a vertex in the embedded and map views}
\label{fig:lensingcube}
\end{figure}

\begin{figure}[htbp]
\centering
\begin{tabular}{c}
\includegraphics[width=0.48\textwidth]
{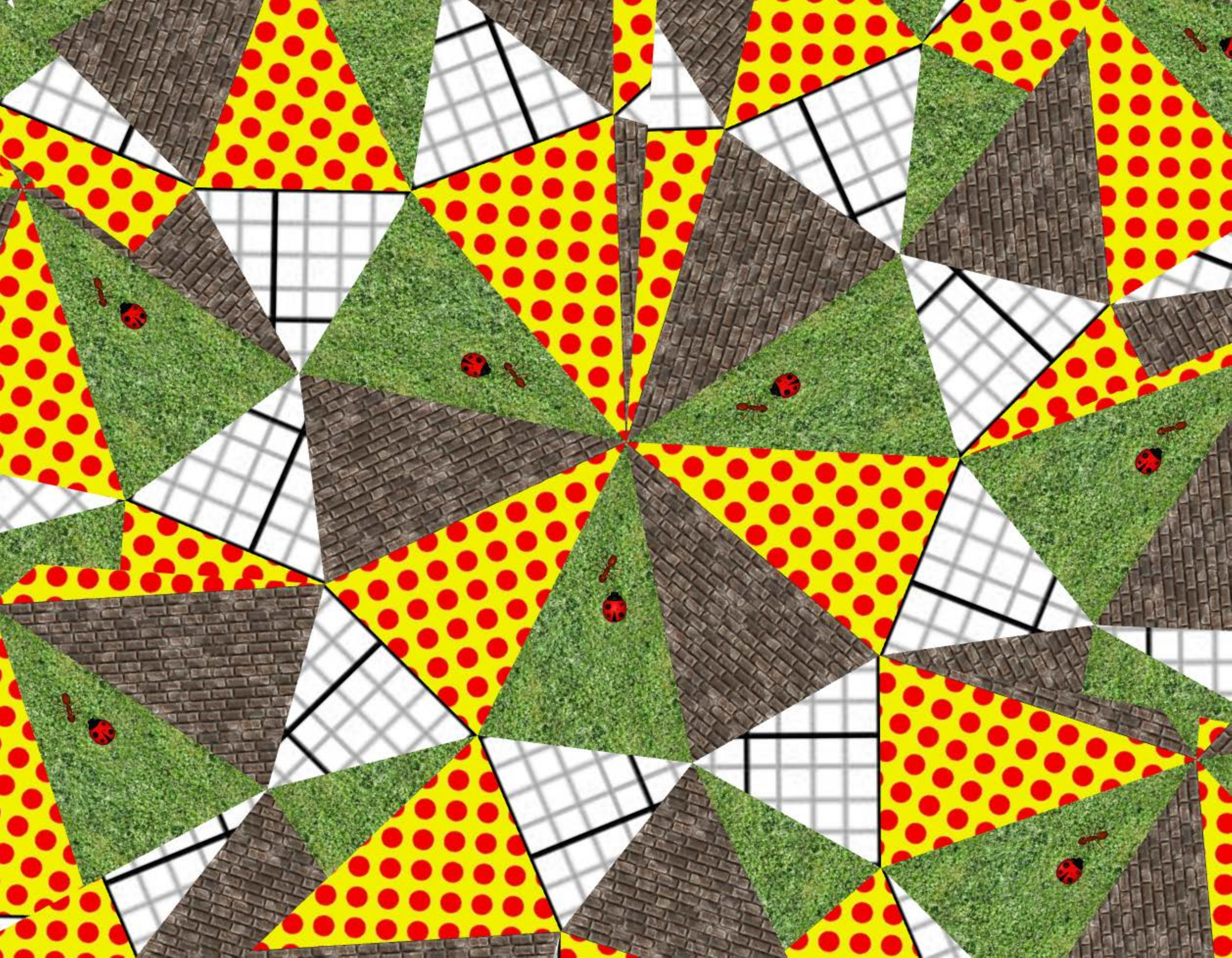}%
\\
\includegraphics[width=0.48\textwidth]
{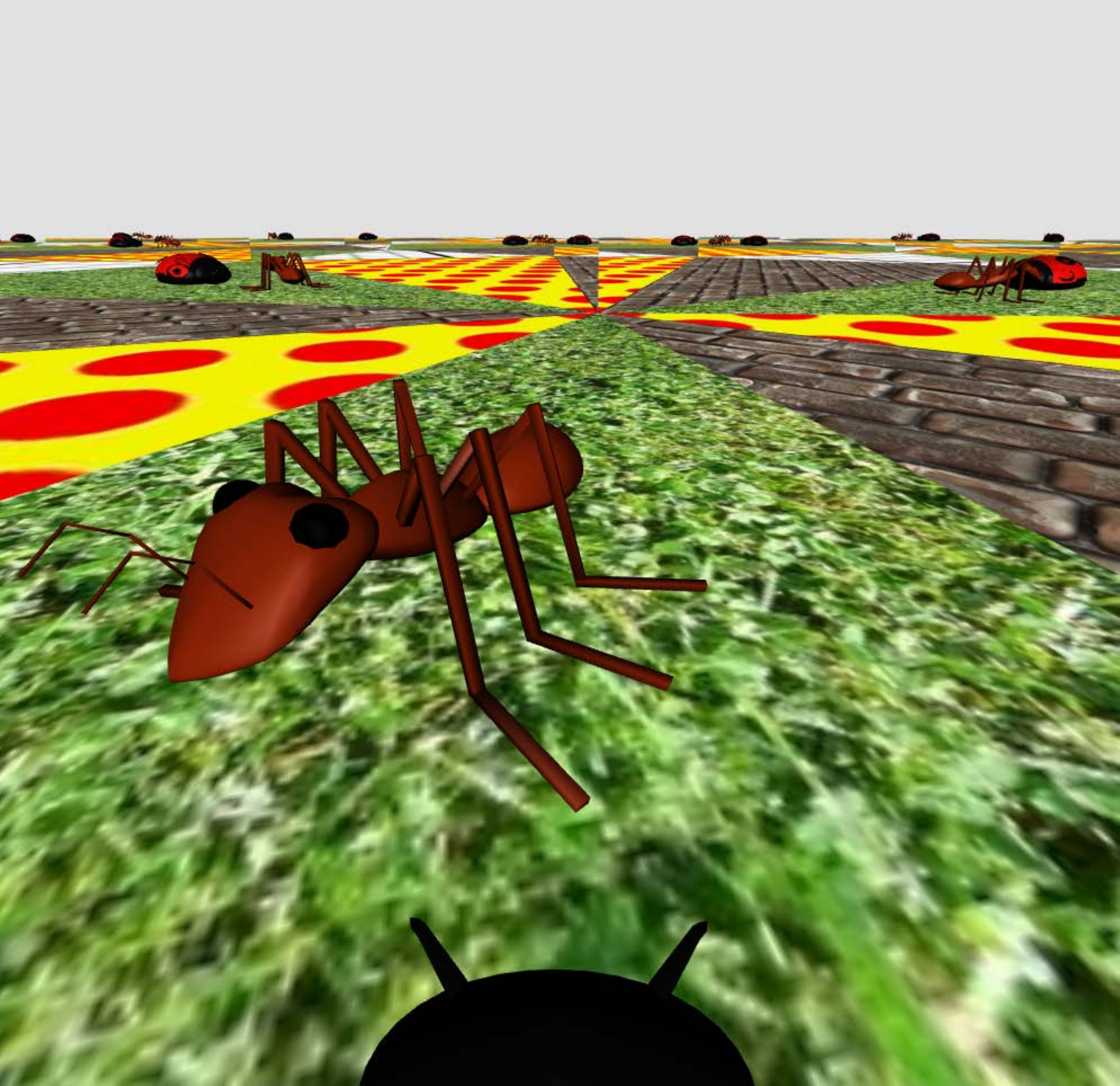}
\end{tabular}
\caption{Lensing on an irregular tetrahedron in map and first person views}
\label{fig:lensingmapfirsttetra}
\end{figure}

A vertex has positive curvature if the sum of the angles of the polygons at
the vertex is less than $2 \pi$.
If one looks in the vicinity of a point with positive curvature, an object on
the other side will appear in more than one direction. Figure \ref{fig:lensingcube} shows how two geodesics can curve around a vertex in a cube to reach the same point. 
If the curvature is large enough, there can be many images
as geodesics can curl around the positive curvature point in multiple ways.
Figure \ref{fig:lensingmapfirsttetra} shows how lensing 
can cause the ladybug to appear to be facing itself. One can also see the lensing of an ant appearing in two places.

\subsection{Cloaking/negative curvature}

\begin{figure}[htbp]
\centering
\begin{tabular}{c}
\includegraphics[width=0.45\textwidth]
{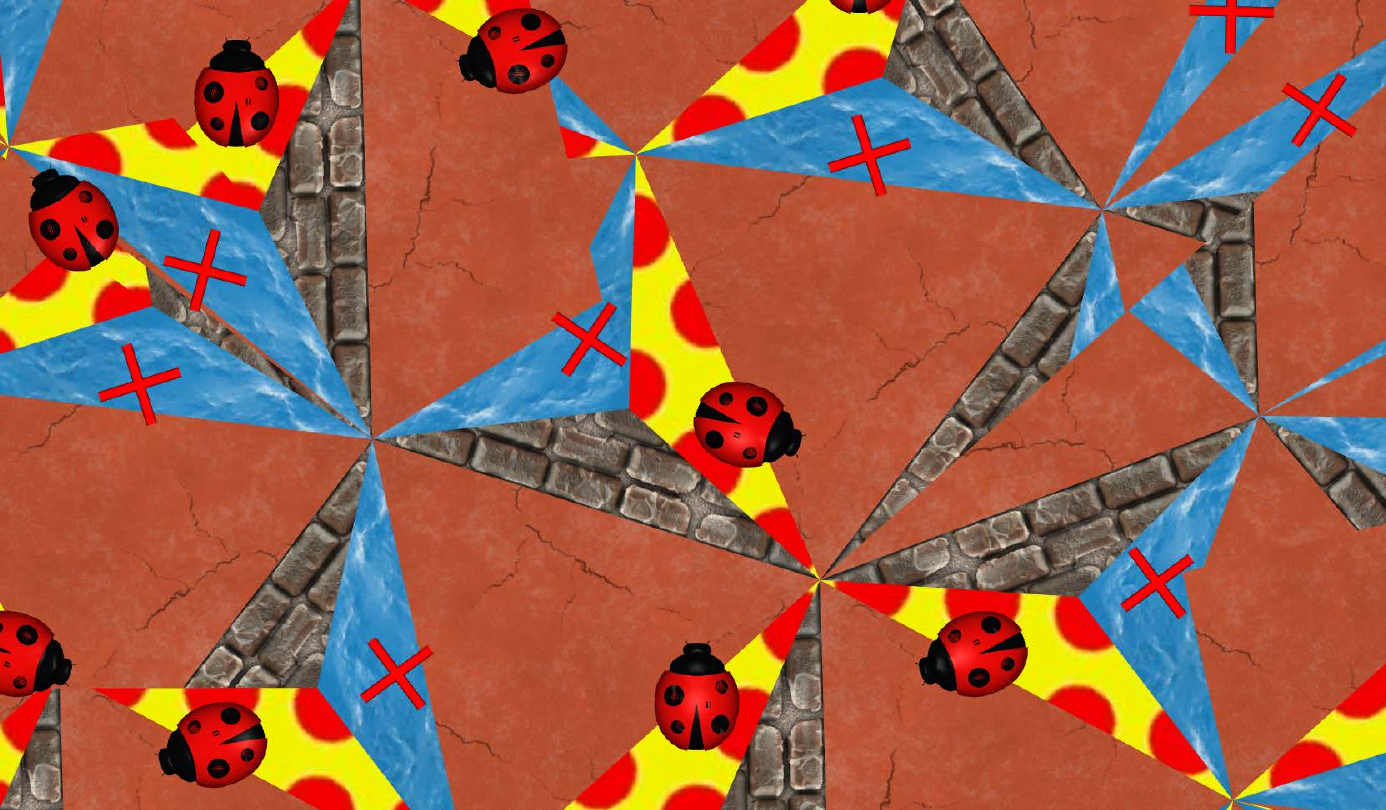}\\\
\includegraphics[width=0.45\textwidth]
{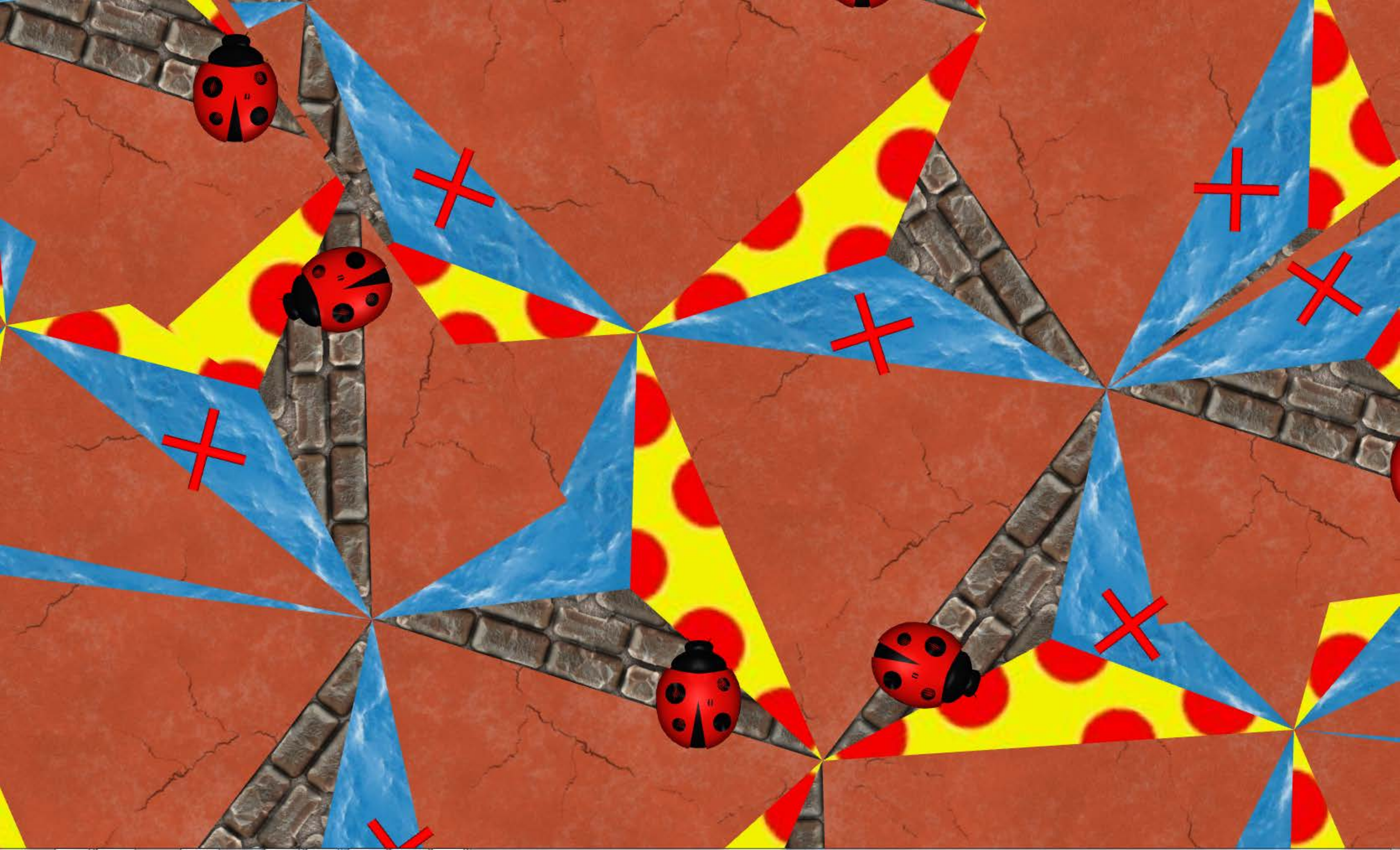}%
\end{tabular}
\caption{Cloaking on map view on a nonembedded tetrahedron}
\label{fig:cloakingmap}
\end{figure}

\begin{figure}[htbp]
\centering
\begin{tabular}{cc}
\includegraphics[width=0.23\textwidth]
{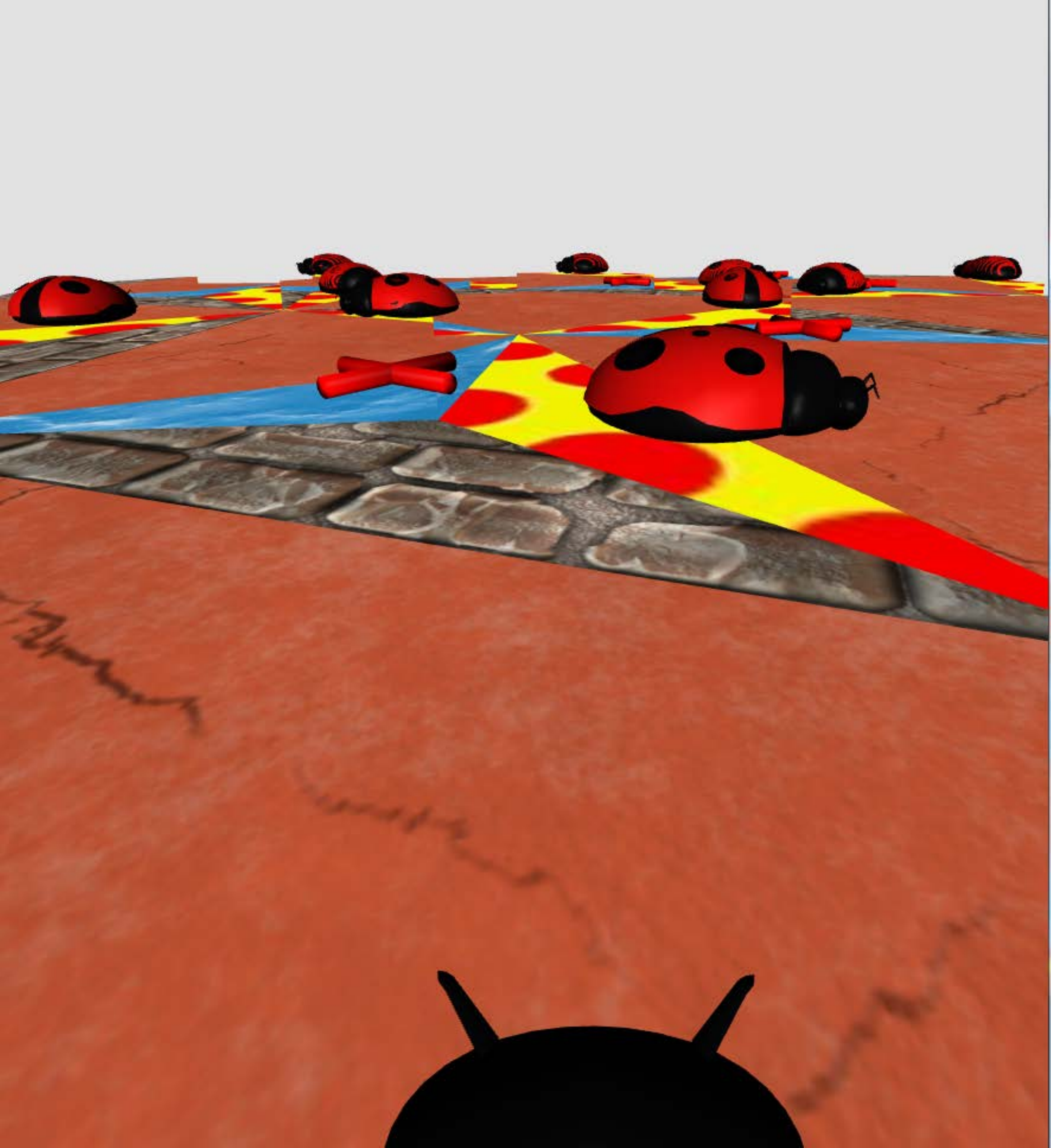}
&
\includegraphics[width=0.23\textwidth]
{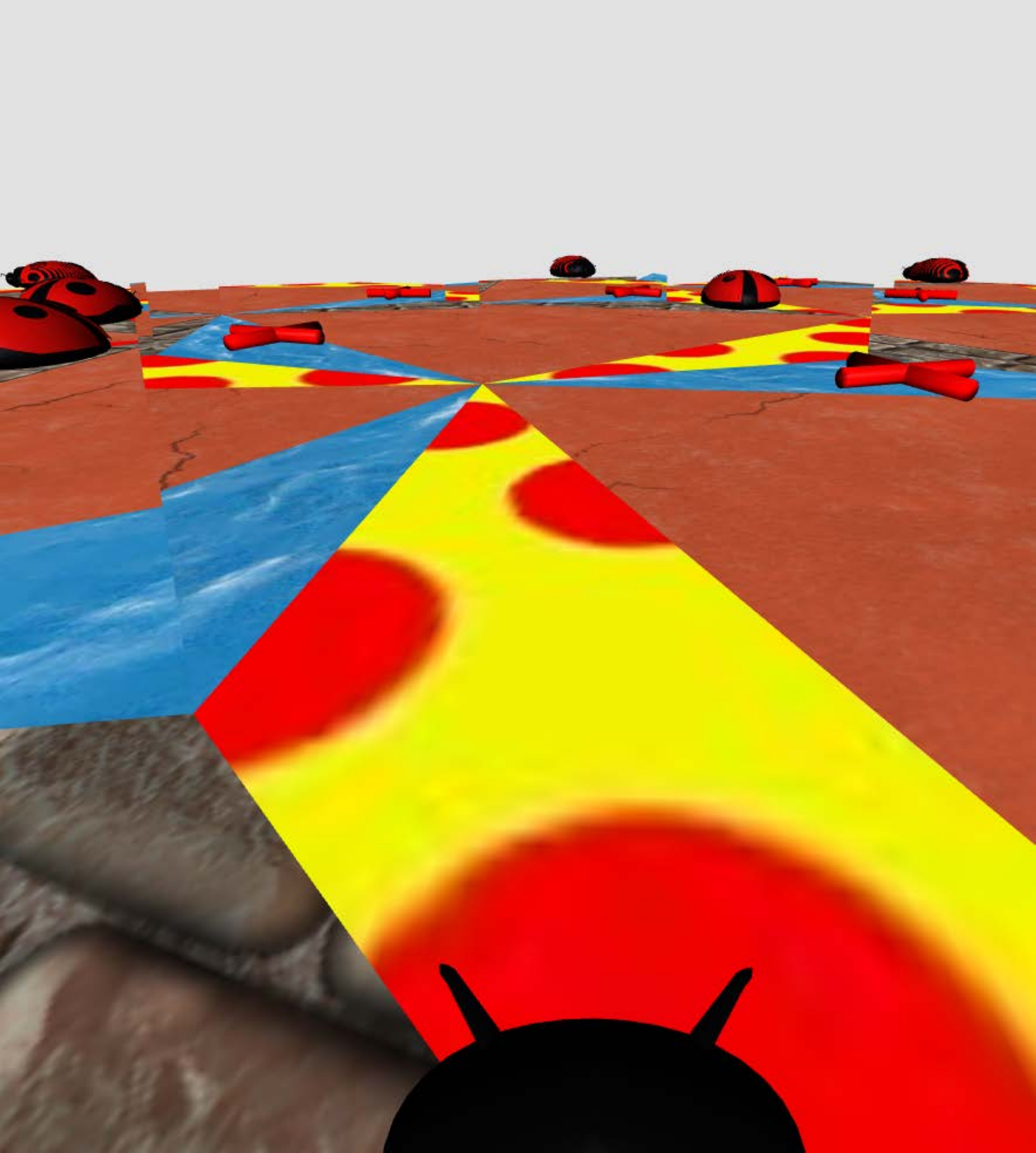}%
\end{tabular}
\caption{Cloaking on first person view on a nonembedded tetrahedron}
\label{fig:cloakingfirstperson}
\end{figure}

A vertex has negative curvature if the sum of the angles of the polygons
at the vertex is more than $2 \pi$. 
In the presence of a point of negative curvature, often all the shortest paths
to a whole region will go through a single point of negative curvature. One
can see this by developing into the plane; the leftover \textquotedblleft
flap\textquotedblright\ of surface will all be accessible through the point of
negative curvature. This means that all light from this region to the viewer
goes through the same point, which essentially means that nothing can be
distinguished from each other on that path, causing that region to be cloaked. 
In order to see objects in that region, one needs
to circle around the point of negative curvature to reveal them. Figures
\ref{fig:cloakingmap} and \ref{fig:cloakingfirstperson} show cloaking; the `X' 
on the blue triangle disappears from view as the ladybug moves forward. This
example is on a tetrahedron made from four triangles that do not fit together
in three-dimensional space since the triangles meeting at one vertex have angles
summing to more than $2 \pi$.
 
Cloaking is due to the fact that there are many geodesics going through a point of 
negative curvature on the polyhedron, and so we cannot distinguish images
on the other side of that point. Note that although one cannot have actual cloaking
 on a smooth manifold, a slightly smoothed surface with lots of
negative curvature in a small region will force all geodesics to go through 
that region. If the region is small enough so that the viewer cannot
distinguish points in the region very well, the same effect will be
achieved.\bigskip

\subsection{Fracturing/nonfracturing}

\begin{figure}[htbp]
\centering
{\includegraphics[width=0.48\textwidth]
{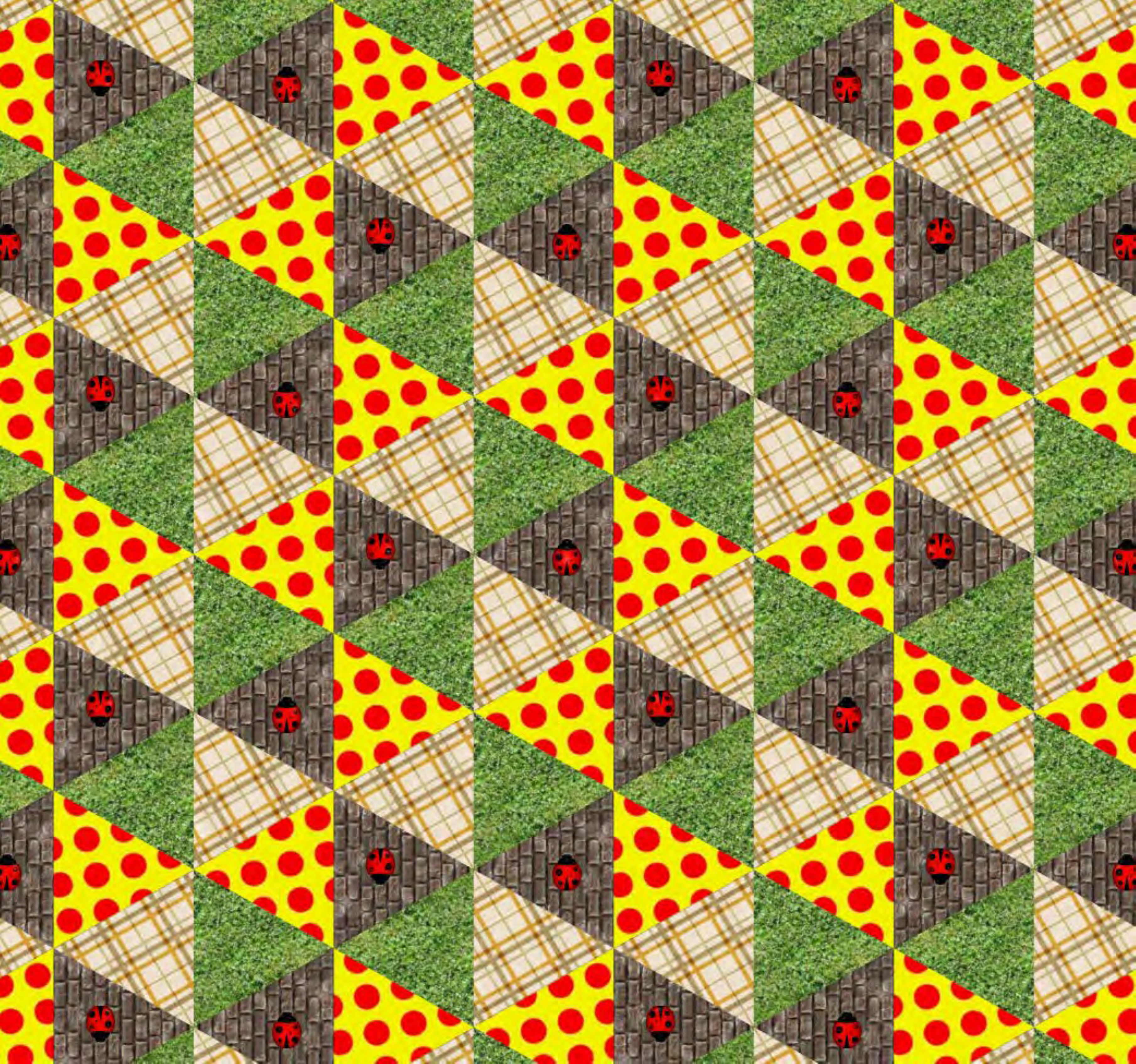}}%
\caption{Non-fracturing of a regular tetrahedron}
\label{fig:nonfractetra}
\end{figure}

\begin{figure}[htbp]
\centering

{\includegraphics[width=0.48\textwidth]
{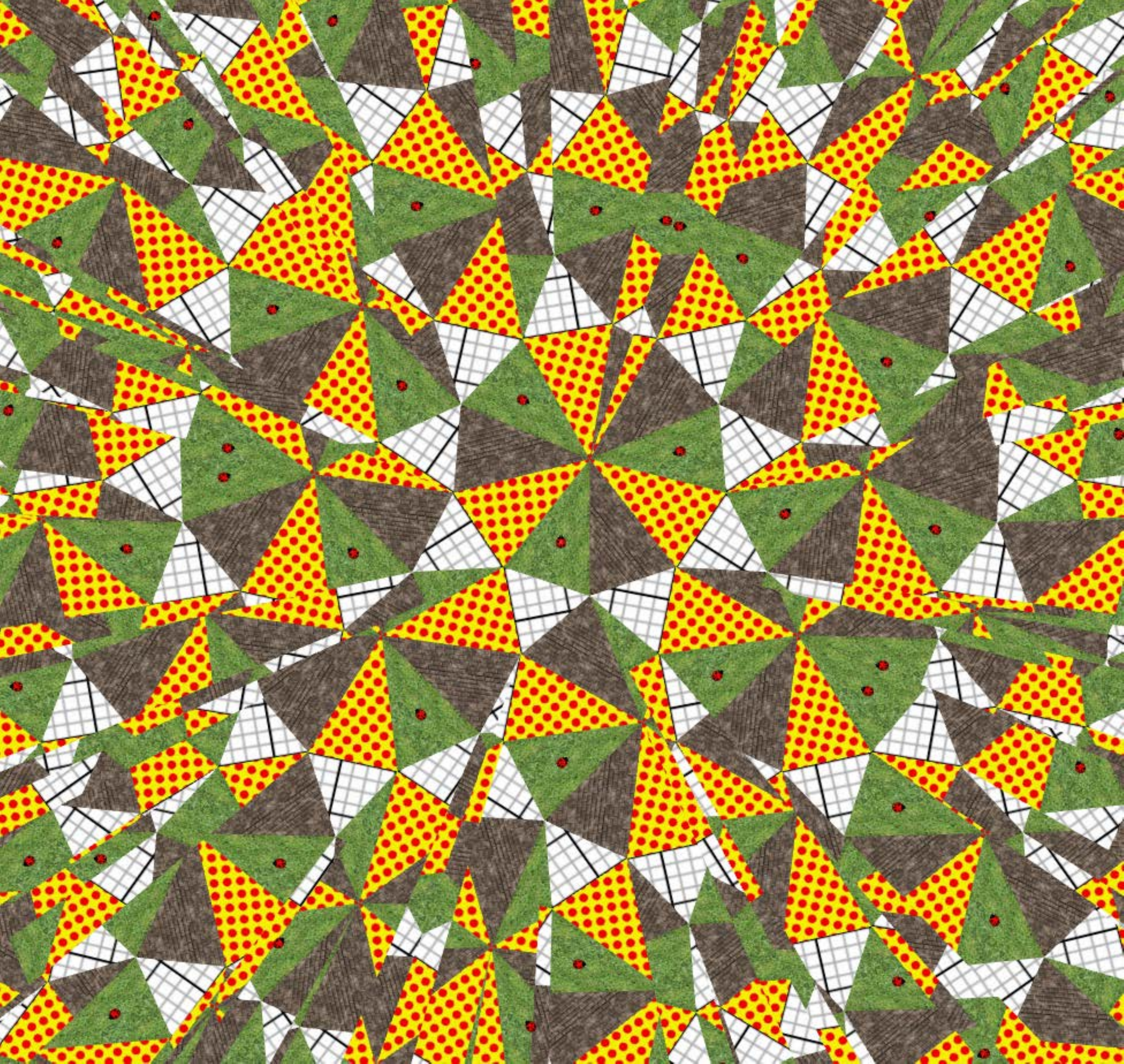}}
\caption{Fracturing of a non-regular tetrahedron}
\label{fig:fractnonregtetra}
\end{figure}

If the sum of the angles at a vertex equal to $2\pi$ divided by an integer, the map
view will not appear to change as one walks through the surface (though there
could be repetition and/or lensing. If it is a rational multiple of $2\pi,$
then the map view will look fairly contiguous, although masking may also
occur. If, however, the sum of angles is not a rational multiple of $2\pi,$
then the as one moves through the manifold, the view will change significantly
along rays that go through vertex points that are not flat. These accumulate
as one goes out a large distance and there appears to be a fracturing of the
manifold in the horizon. Figure \ref{fig:nonfractetra} shows the non-fracturing case of the
regular tetrahedron and Figure \ref{fig:fractnonregtetra} shows fracturing.\bigskip


\subsection{Apexing}

\begin{figure}[htbp]
\begin{tabular}{ccc}
\includegraphics[width=0.15\textwidth]
{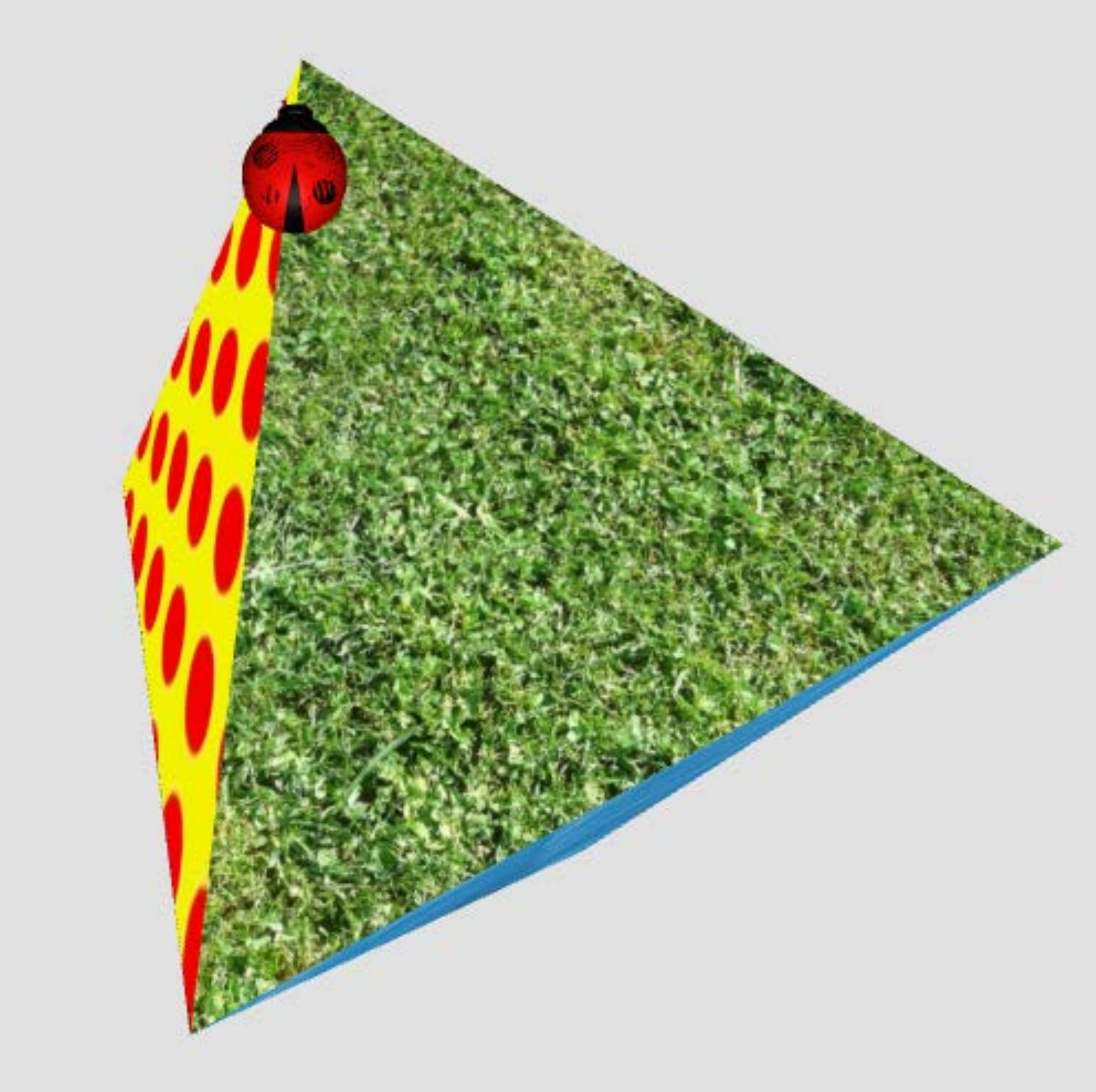}
&
\includegraphics[width=0.15\textwidth]
{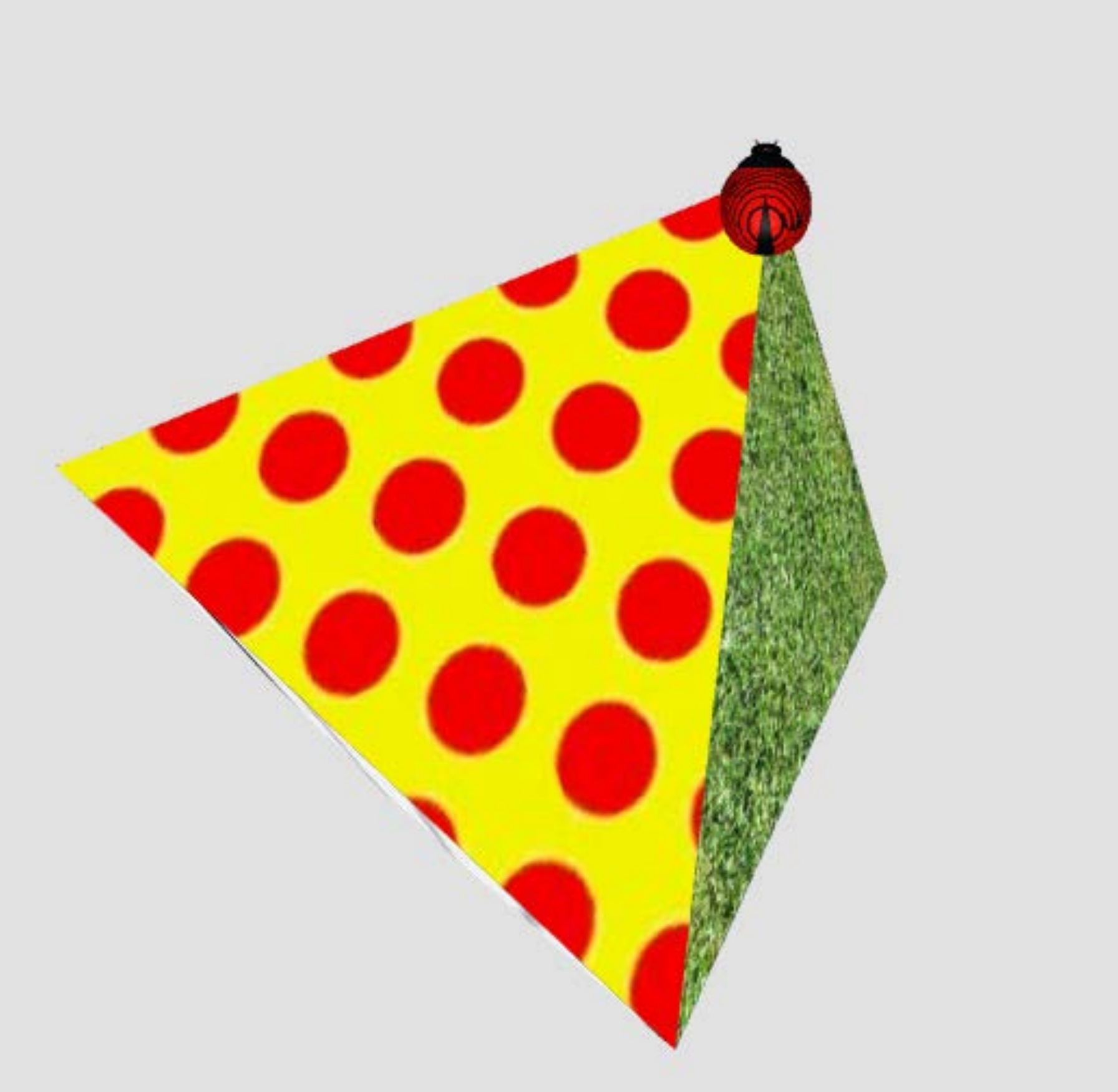}
&
\includegraphics[width=0.15\textwidth]
{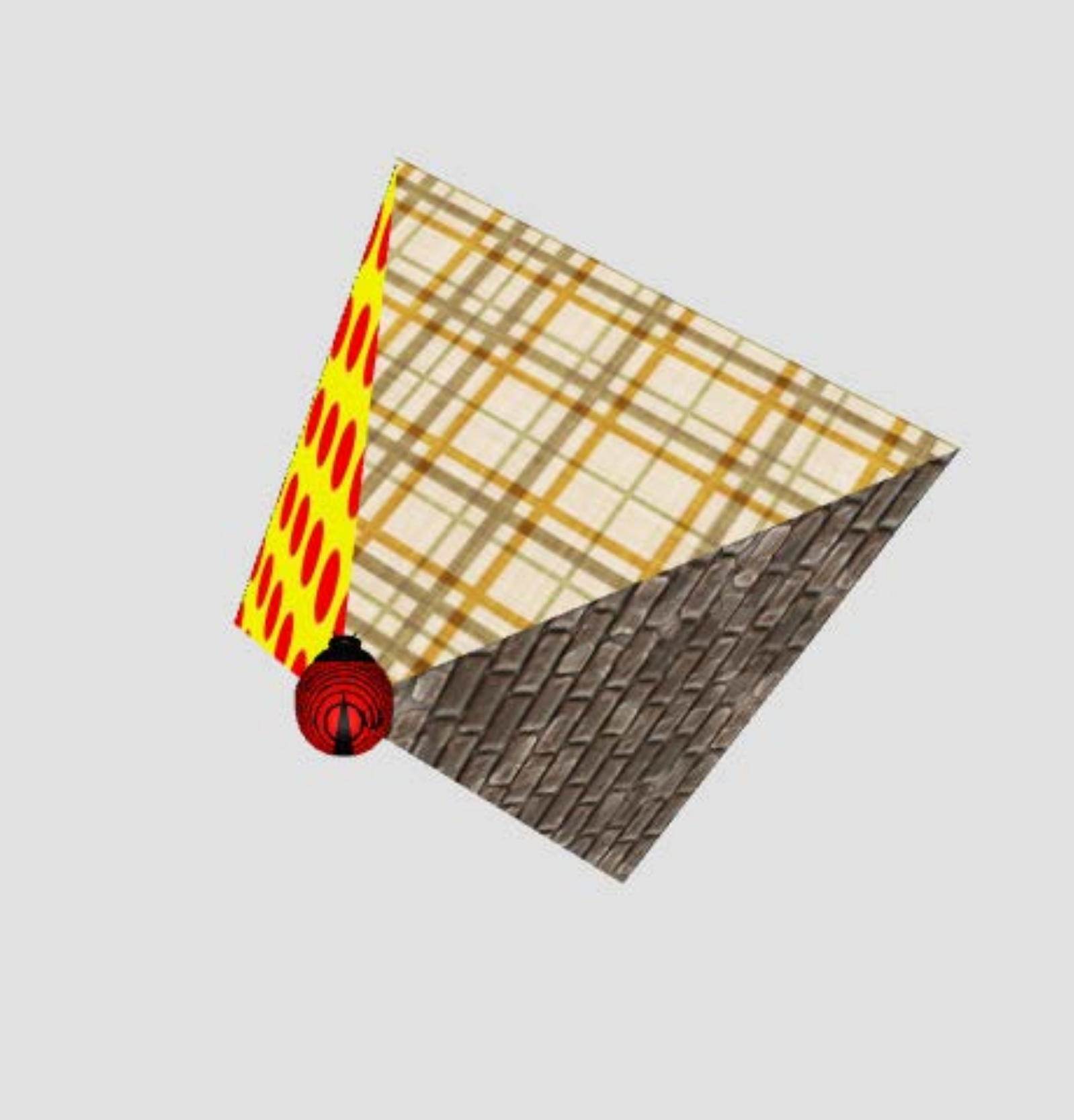}%
\\
\includegraphics[width=0.15\textwidth]
{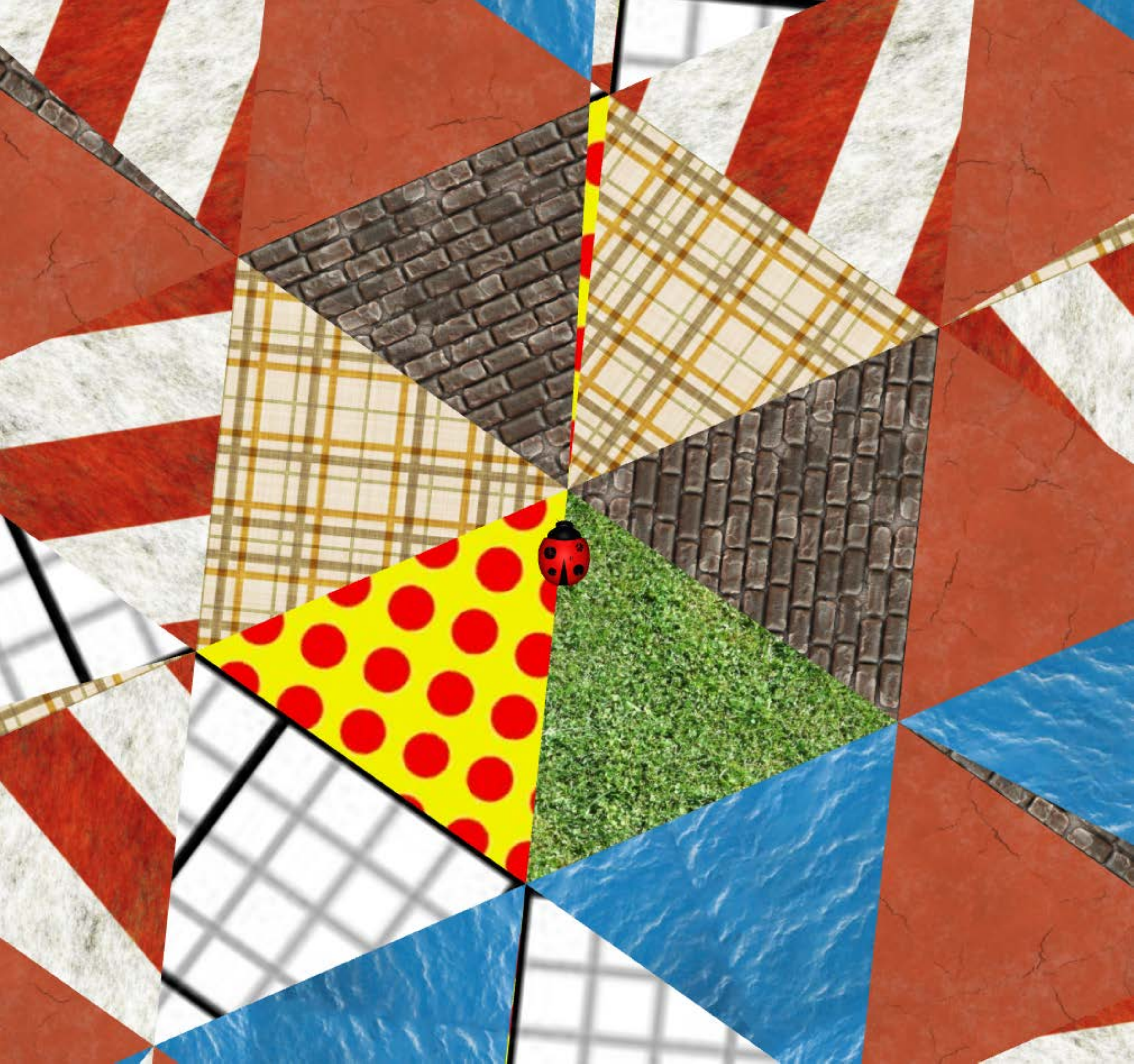}
&
\includegraphics[width=0.15\textwidth]
{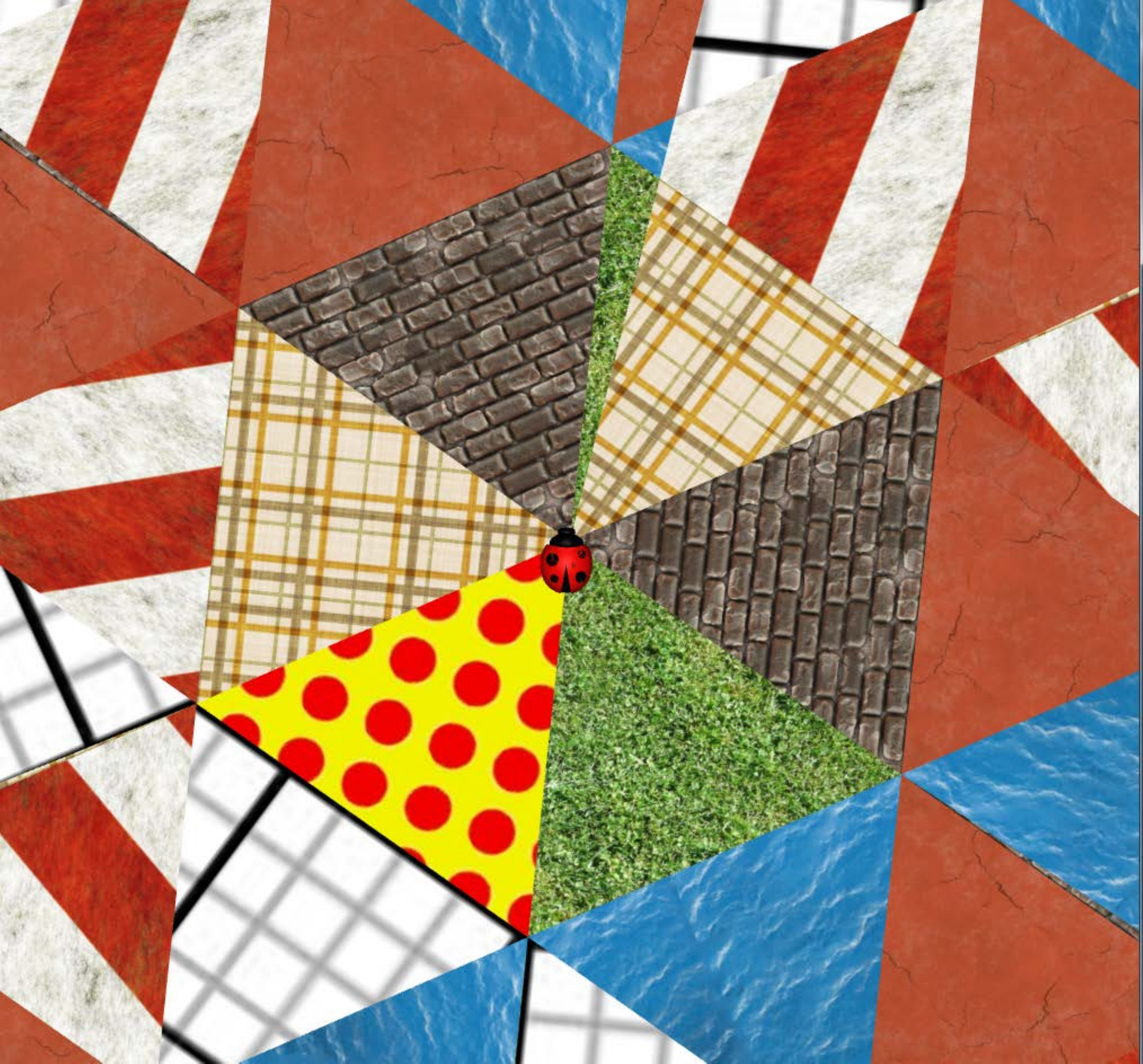}
&
\includegraphics[width=0.15\textwidth]
{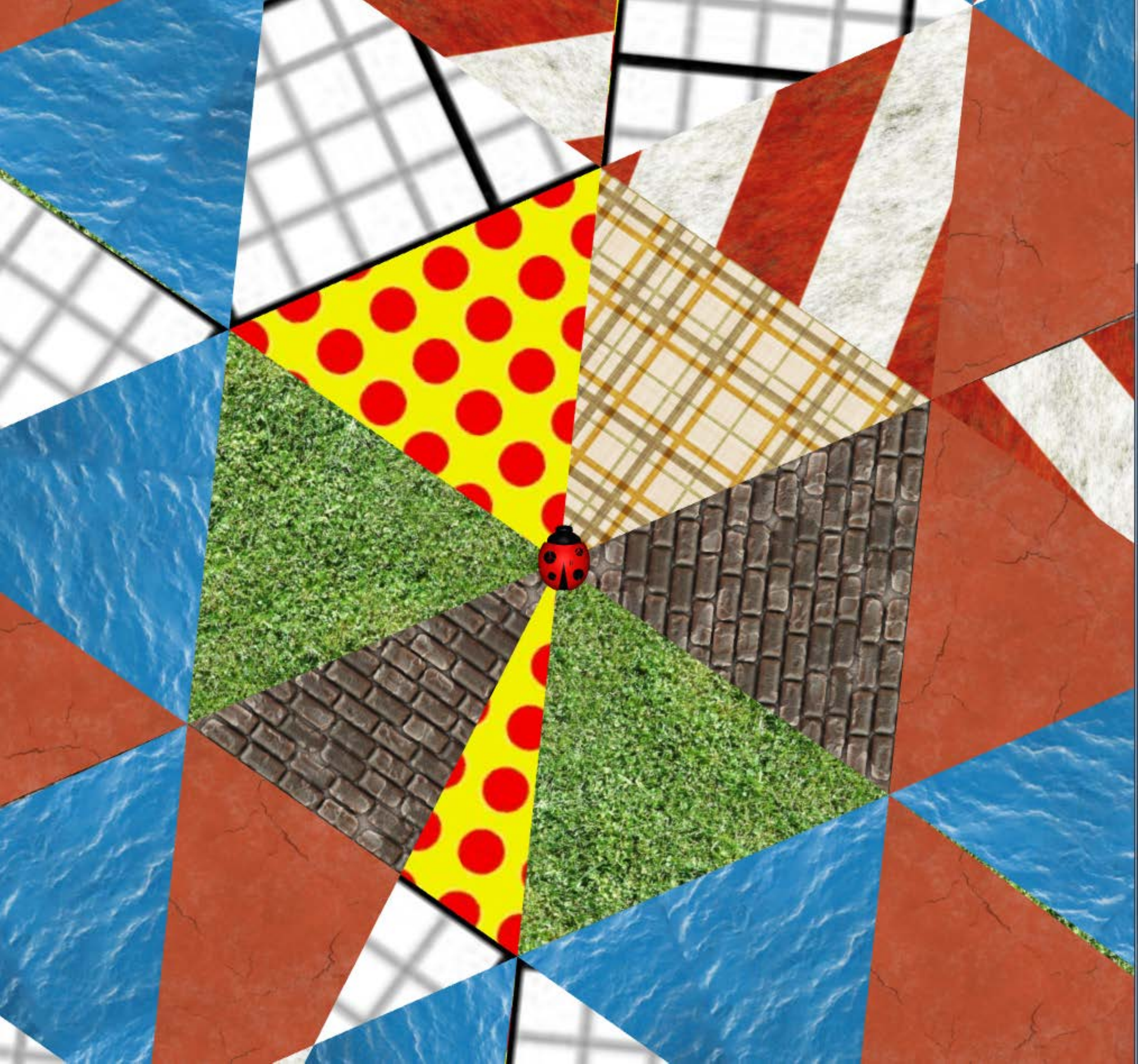}%
\\
\includegraphics[width=0.15\textwidth]
{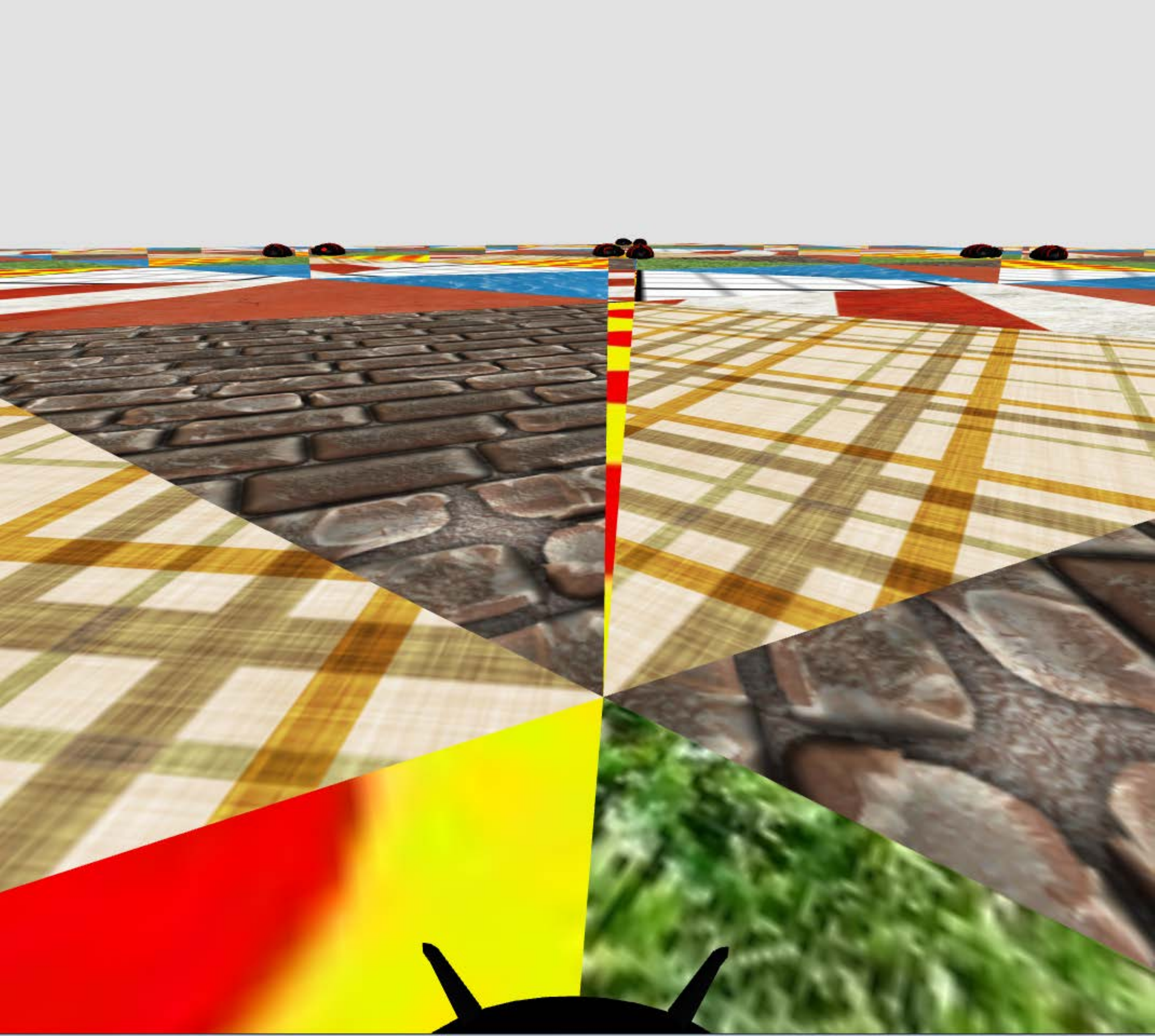}
&
\includegraphics[width=0.15\textwidth]
{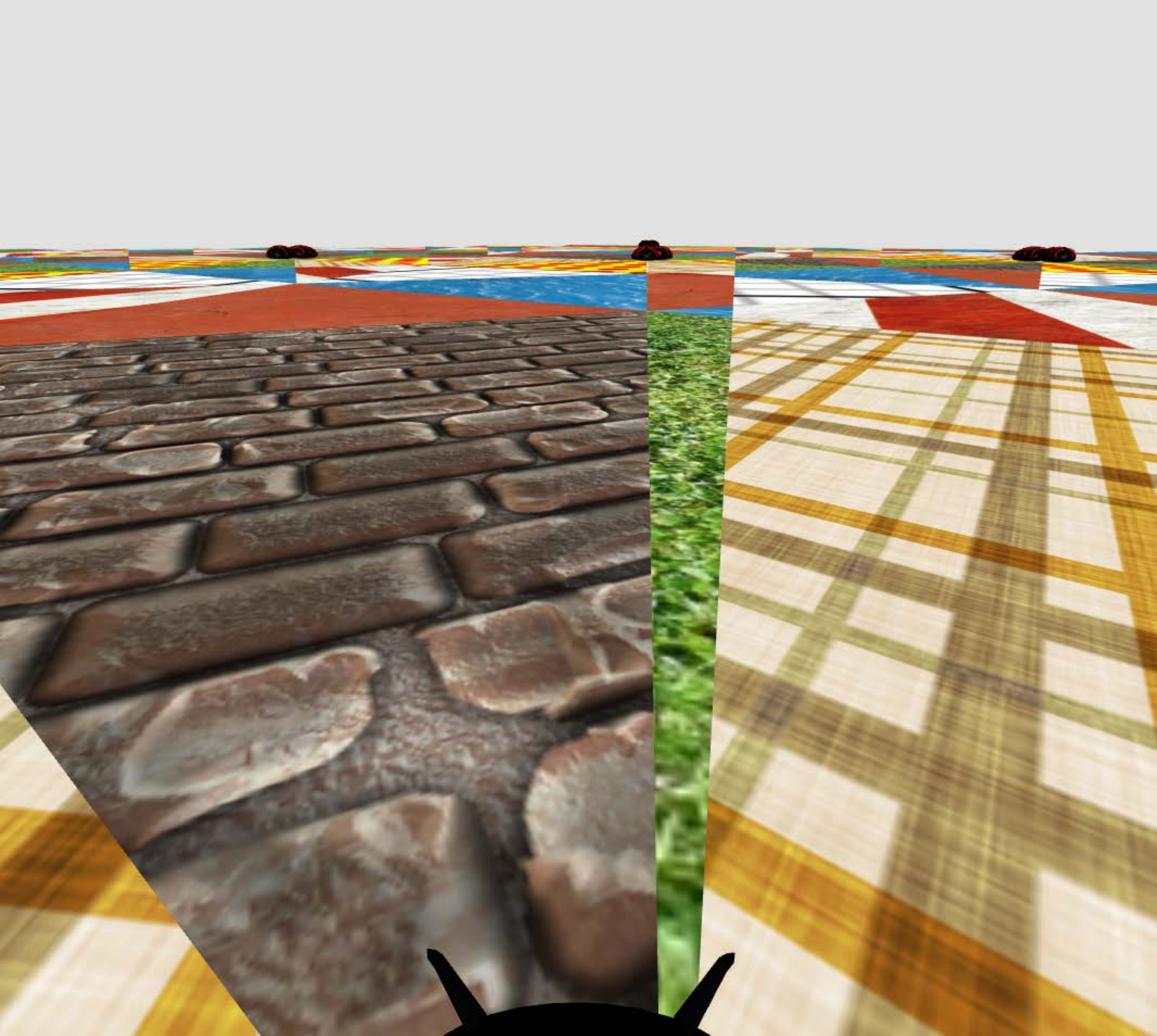}%
&
\includegraphics[width=0.15\textwidth]
{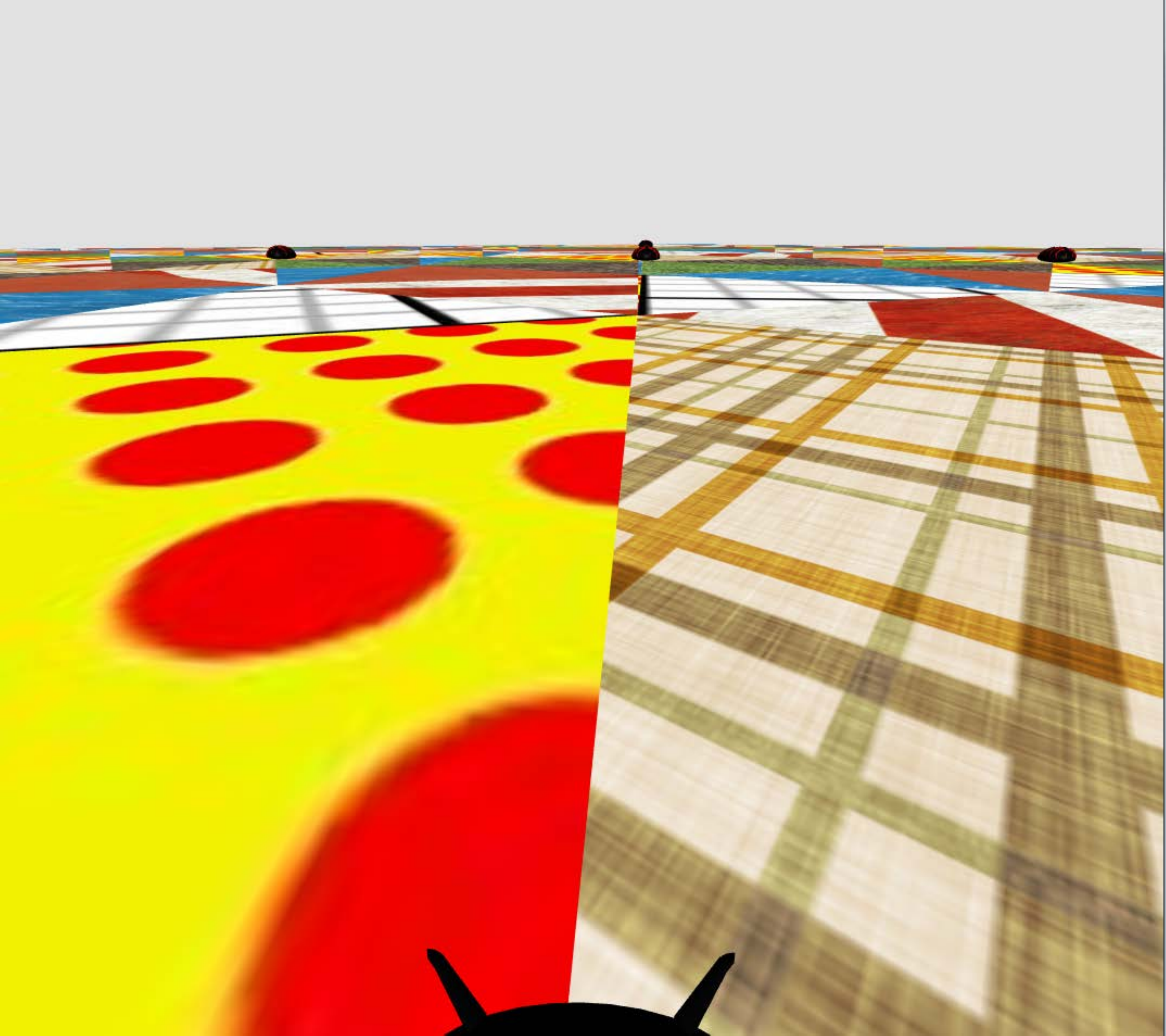}%
\end{tabular}
\caption{Apexing in embedded, map, and first person views on an octahedron}
\label{fig:apex}
\end{figure}

Since there are no shortest paths that go through a vertex of positive
curvature and there are many shortest paths that go through a vertex of
negative curvature, when walking through a vertex the resulting position and
orientation can be a little unpredictable. This is analogous to trying to run
over the apex of a mountain; it is not always clear how you will be oriented
when you go down. See Figure \ref{fig:apex}. 

\section{Conclusion}
We have seen several way that the geometry of polyhedra affect how light travels 
around a surface that are quite different than what happens in the Euclidean plane.
 Finite extent can cause straight lines to curve back on themselves, resulting in 
 repetition of images. Positive curvature can cause straight lines to swirl around 
 to avoid the positive curvature point, which we call lensing. Negative curvature 
 can cause many straight lines to go through the same negative curvature point, and 
 objects on the other side of that point can be cloaked. Finally, the fact that there 
 is curvature at all can cause the size of the pieces of a polygon that appear together
  to be fractured at large distances unless the curvatures are especially nice (for 
  instance, a rational multiple of $2 \pi$) and causes unpredictable behavior
  when moving through a vertex with curvature.

Many of these phenomena occur on smooth Riemannian manifolds due to finite extent in some direction
(short geodesics) and curvature. The main difference with polyhedra is that the curvature
is concentrated at vertex points, where there is infinite curvature at isolated points.

\section*{Acknowledgements}
The software used to generate the pictures in this article allows for interactive exploration
of polyhedral surfaces, both embedded and not. It was developed by a number
of graduate and undergraduate students, most notably graduate students Thomas (Danny) Maienschein
and Joseph Thomas and undergraduate researchers Joseph Crouch, Mark Doss, Taylor Johnson, Kira Kiviat, Justin Lanier, Taylor Corcoran, Qiming Shao, Staci Smith, Jeremy Mirchandani, and Tanner Prynn. Much thanks to previous developers and collaborators Dan Champion, Yuliya Gorlina, Alex Henniges, Tom Williams, Mitch Wilson, Kurtis Norwood, and Howard Cheng. The software is available on github \cite{github} and uses jReality. This work was supported by NSF DMS 0748283 and DMS 0602173.

%
%
%
%
%

\bibliographystyle{plain}

\end{document}